\documentclass[12pt]{amsart}
\usepackage{amscd, amsmath, amstext, amssymb, amsbsy}
\theoremstyle{plain}
\newtheorem{thm}{Theorem}[section]
\newtheorem{lemma}[thm]{Lemma}
\newtheorem{prop}[thm]{Proposition}
\newtheorem{question}[thm]{Question}
\newtheorem{conj}[thm]{Conjecture}
\theoremstyle{definition}
\newtheorem{exam}[thm]{Example}

\newcommand{\ol}{\overline}
\newcommand{\s}[1]{\mathcal{#1}}
\newcommand{\appa}{\subseteq}
\newcommand{\BB}{\mathcal{B}}
\newcommand{\DD}{\mathcal{D}}
\newcommand{\EE}{\mathcal{E}}
\newcommand{\MM}{\mathcal{M}}
\newcommand{\cxd}{\mathbb{C}^{d}}
\newcommand{\mbb}{\mathbb}
\newcommand{\la}{\left<}
\newcommand{\ra}{\right>}
\newcommand{\rank}{\mathrm{rank} \, }
\newcommand{\trace}{\mathrm{trace} \, }
\newcommand{\ran}{\mathrm{ran} \, }
\newcommand{\hd}{H_{d}^{2}}
\newcommand{\HH}{\mathcal{H}}
\newcommand{\lra}{\longrightarrow}

\newcommand{\npker}[2]{1-\frac{1}{k_{#1}(#2)}}

\newcommand{\sot}{\; \text{(SOT)} \; }
\newcommand{\sra}{\rightarrow}
\renewcommand{\Re}{\text{Re} \, }
\newcommand{\prtl}{\partial}
\newcommand{\ds}{\displaystyle}

\newcommand{\wot}{\; \text{(WOT)} \; }
\begin{document}
\bibliographystyle{amsplain}

\title[The Structure of Inner Multipliers]{The 
Structure of Inner Multipliers on Spaces with Complete Nevanlinna Pick 
Kernels}
%\author{Devin Greene, Stefan Richter, and Carl Sundberg}
\author{Devin~C.V. Greene}
\address{Department of Mathematics\\
University of Calfornia\\
Berkeley, CA 94720}
\email{deving@math.berkeley.edu}
\author{Stefan Richter}
\address{Department of Mathematics\\
University of Tennessee\\
Knoxville, TN 37996-1300}
\email{Richter@math.utk.edu}
\author{Carl Sundberg}
\email{Sundberg@math.utk.edu}
\thanks{Work of the
second and third author was supported by the National
Science Foundation}

\begin{abstract}
Let $k$ be the reporducing kernel for a Hilbert space $\mathcal{H}(k)$
of nanlytic functions on $B_{d}$, the open unit ball in
$\mathbb{C}^{d}$, $d \geq 1$.  $k$ is called a complete NP kernel, if
$k_{0} \equiv 1$ and if $1-1/k_{\lambda}(z)$ is positive definite on
$B_{d} \times B_{d}$.  Let $\mathcal{D}$ be a separable Hilbert space,
and consider $\HH(k) \otimes \DD \cong \HH(k, \DD)$, and think of it
as a space of $\DD$-valued $\HH(k)$-functions.  A theorem of
McCullough and Trent, \cite{MT}, partially extends the
Beurling-Lax-Halmos theorem for the invariant subspaces of the Hardy
space $H^{2}(\DD)$.  They show that if $k$ is a complete NP kernel and
if $\DD$ is a separable Hilbert space, then for any scalar multiplier
invariant subspace $\MM$ of $\HH(k, \DD)$ there exists an auxiliary
Hilbert space $\EE$ and as multiplication operator $\Phi: \HH(k , \EE)
\lra \HH(k, \DD)$ such that $\Phi$ is a partial isometry and $\MM =
\Phi \HH(k,\EE)$.  Such multiplication operators are called inner
multiplication operators and they satisfy $\Phi \Phi^{*} = $ the
projection onto $\MM$.

In this paper we shall show that for many interesting complete NP kernels the analogy 
with the Beurling-Lax-Halmos theorem can be strengthened.  We show that for almost 
every $z \in B_{d}$ the nontangential limit $\phi(z)$ of the multiplier $\phi: B_{d} 
\lra \BB(\EE, \DD)$ associated with an inner multiplication operator $\Phi$ is a 
partial isometry and that $\rank \phi(z)$ is equal to a constant almost everywhere.  

The result applies to certain weighted Dirichlet spaces and to the symmetric Fock 
space $H_d^2$.  In particular, our result implies that the curvature 
invariant of W. 
Arveson (\cite{arvcurv}) of a pure contractive Hilbert module of finite rank is an 
integer.  The answers a question of W. Arveson, \cite{arvcurv}.
\end{abstract}

\maketitle

\section{Introduction} \label{1}

For a positive integer $d$ we denote the unit ball in $\cxd$ by $B_{d} = \{\lambda 
\in \cxd: |\lambda| < 1\}$. If $\HH$ is a Hilbert space of analytic functions 
on 
$B_{d}$ such that for each $\lambda \in B_{d}$ the point evaluation $\lambda 
\mapsto 
f(\lambda)$ is a continuous linear functional on $\HH$, then $\HH$ has a reproducing 
kernel $k$; that is, for each $\lambda \in B_{d}$ there is a $k_{\lambda} \in \HH$ 
such that $f(\lambda) = \la f, k_{\lambda} \ra$ for each $f \in \HH$.  As a function 
of $\lambda$ and $z$ in $B_{d}$, $k_{\lambda}(z)$ is a positive definite function 
which is analytic in $z$ and anti-analytic in $\lambda$.  It is well known that $k$ 
determines the space $\HH$.  Thus, we shall write $\HH(k)$ for the space of analytic 
functions with reproducing kernel $k$.  

An analytic function $\varphi$ on $B_{d}$ is a multiplier of $\HH(k)$ if $\varphi f 
\in \HH(k)$ for every $f \in \HH(k)$.  We shall write $M(k)$ for the collection of 
all multipliers.  A standard argument with the closed graph theorem shows that each 
$\varphi \in M(k)$ defines a bounded linear operator $M_{\varphi}: f \mapsto 
\varphi f$ on $\HH(k)$.  Thus we define the multipier norm by $\|\varphi\|_{M} = 
\|M_{\varphi}\|$.  A subspace $\MM$ of $\HH(k)$ is called a multiplier invariant 
subspace if $\varphi \MM \subseteq \MM$, and we shall denote the collection of all 
multiplier invariant subspaces by $\mathrm{Lat} \, M(k)$.  

A reproducing kernel $k$ on $B_{d}$ is called a {\it complete Nevanlinna-Pick kernel} 
({\it complete NP kernel} for short), if $k_{0}(z) = 1$ for all $z \in B_{d}$ and if 
there exists a sequence of analytic functions $\{b_{n}\}_{n \geq 1}$ on $B_{d}$ such 
that 

\begin{equation} \label{1.1}
\npker{\lambda}{z} = \sum_{n \geq 1} \, b_{n}(z) \ol{b_{n}(\lambda)} 
\text{~for all~} \lambda, \, z \in B_{d}.
\end{equation}

We note that the sequence may be finite, and that this condition is actually 
equivalent to the assumption that $1-1/k$ is positive definite.  Complete NP kernels 
have been investigated in connection with Nevanlinna-Pick and Caratheodory 
interpolation and commutant lifting properties (see \cite{agler2}, 
\cite{aglermccarthy1}, \cite{balltrentvinnikov}, \cite{mccullough1}, 
\cite{quiggin1}).  Examples of spaces 
with complete NP kernels on the unit disc $\mbb{D} = B_1$ are the Hardy space $H^2 
(\mbb{D})$ of the unit disc $\mbb{D}$, the Dirichlet space of all analytic 
functions 
on $\mbb{D}$ whose derivative is square area integrable, or, more generally, the 
weighted Dirichlet spaces $D_{\alpha}, \, \alpha \geq 0$ on the unit disc (see 
Example~\ref{Example 4.4} for definitions and details).  For $d \geq 1$ we 
mention the space $\hd$ on $B_d$, 
which is defined by the kernel $k_\lambda(z) = \frac{1}{1-\la z, \lambda \ra_d}$, 
$\la 
z,\lambda \ra_d = \sum_{i=1}^d \, z_i \ol{\lambda_i}$.  The space $\hd$ was 
investigated in \cite{arvmvot}, \cite{arvcurvinv}, and \cite{arvcurv}, because of its 
connection to the dilation theory of certain commuting operator tuples, so called 
$d$-contractions, or row contractions.  In Proposition~\ref{Proposition 2.3}
we shall see that for all complete NP kernels $k$ one has $\HH(k) \subseteq 
H^2(\partial B_d)$, the ordinary Hardy space of the unit ball.  But we note that for 
$d > 1$, the reproducing kernel for $H^2(\partial B_d)$ is not a complete NP kernel.  

We shall now fix a complete NP kernel $k$ and a sequence $\{b_n\}$ as in 
\eqref{1.1}.  
One shows that each $b_n \in M(k)$ and that $P_0 = \sum_{n \geq 1} \, M_{b_n} 
M_{b_n}^*$ $\sot$, where $P_0$ is the projection onto the multiplier 
invariant subspace 
$\HH_0 = \{f \in \HH(k):f(0) = 0\}$ (see Lemma 1.4 of \cite{MT}).  It 
is remarkable that it follows that the projection onto every multiplier invariant 
subspace can be written in a similar manner.  the general case of the 
following 
theorem is due to McCullough and Trent, \cite{MT}, and for the special 
case 
of $\HH(k) = \hd$ it was found by Arveson in \cite{arvcurvinv}.  

\begin{thm} \label{Theorem 1.1}
Let $k$ be a complete NP kernel and let $\MM$ be a multiplier invariant subspace.  
Then there exists a sequence $\{\varphi_n\} \subseteq M(k) \cap \MM$ such that 

\begin{equation} \label{1.2}
P_\MM = \sum_{n\geq1} \, M_{\varphi_n} M_{\varphi_n}^* \sot,
\end{equation}

\noindent where $P_\MM$ is the projection onto $\MM$.

\end{thm}

We make several remarks.  First, McCullough and Trent prove this theorem in 
a 
somewhat more general setting; it is not even necessary that the kernel $k$ is 
reporducing for a space of {\it analytic} functions.  Secondly, by applying the 
expression \eqref{1.2} to the reproducing kernel $k_\lambda, \, \lambda \in B_d$, 
one 
obtains 

\begin{equation} \label{1.3}
\sum_{n \geq 1} \, |\varphi_n(\lambda)|^2 = \frac{\|P_\MM 
k_\lambda\|^2}{\|k_\lambda\|^2} \leq 1.
\end{equation}

Thus, each function $\varphi_n$ is in the unit ball of $H^\infty(B_d)$ and therefore 
for a.e. $z \in \partial B_d$ the nontangential limit $\varphi_n(z)$ of $\varphi$ 
exists.  Here, and in what follows a.e. stands for a.e. $[\sigma]$ where $\sigma$ is 
the rotationally invariant probability measure on $\partial B_d$.  

Of course, Beurling's theorem implies that for $H^2(\mbb{D})$ the sequence 
$\{\varphi_n\}$ can be chosen to be a single inner function $\varphi$, which 
satisfies $|\varphi(z)| = 1$ for a.e. $z \in \partial \mbb{D}$.  Our first main 
result is the following.

\begin{thm} \label{Theorem 1.2}
Let $k$ be a complete NP kernel on $B_d$ and assume that there is a set $\s{P} 
\appa \s{H}(k) \cap C(\ol{B_d})$ which is dense in $\s{H}(k)$ and such that for all 
$p 
\in \s{P}$ and $z \in \partial B_d$, $\lim_{\lambda \rightarrow z} 
\frac{\|pk_\lambda\|}{\|k_\lambda\|} = |p(z)|$.  Then any sequence $\{\varphi_n\}$ 
which is associated with a nonzero multiplier invariant subspace $\s{M}$ as in 
\eqref{Theorem 1.1} is an inner sequence, i.e.

\begin{equation*}
\sum_{n \geq 1} \, |\varphi_n (z)|^2 = 1 \quad \text{for} \; \; [\sigma]  
\; \; \text{a.e.} 
\; \; z \in \partial B_d.
\end{equation*}
\end{thm}

We shall show in Section 4 that the hypothesis of the theorem is satisfied for many 
complete NP kernels that are invariant under the unitary group $\s{U}$.  We shall 
call $k$ unitary invariant, or $\s{U}$-invariant, if $k_\lambda (z) = f(\la 
z,\lambda \ra_d)$ for some function $f$ on the unit disc.  One can show that the 
$\s{U}$-invariant kernels are precisely the kernels that satisfy $k_{U\lambda} (Uz) = 
k_\lambda(z)$ for each unitary map $U:\mbb{C}^d \lra \mbb{C}^d$.  More precisely, we 
shall see that whenever a complete NP kernel $k$ is of the form 

\begin{equation*}
k_\lambda(z) = \sum_{n=0}^\infty \, a_n (\la z, \lambda \ra _d )^n,
\end{equation*}

\noindent where $a_n > 0$, $\sum_{n=0}^\infty \, a_n = \infty$, and $\lim_{n 
\sra \infty} \, a_n/a_{n+1} = 1$, then the hypothesis of Theorem~\ref{Theorem 
1.2} is 
satisfied.  In particular, the theorem applies to the weighted Dirichlet spaces 
$D_\alpha$, $0 \leq \alpha \leq 1$, and to the space $H_d^2$, and we shall 
see that 
if $k$ is $\s{U}$-invariant then the multiplier invariant subspaces are exactly the 
subspaces which are invariant under the multiplication by all the coordinate 
functions $z \mapsto z_i$, $i = 1, \ldots, d$.  

However, we shall see that the conclusion of Theorem~\ref{Theorem 1.2} does 
not hold for 
the 
weighted Dirichlet spaces $D_\alpha$, $\alpha > 1$.  

For the space $H_d^2$ this theorem was conjectured by Arveson, 
\cite{arvcurv}, 
\cite{arvcurvinv}.  He proved the theorem for invariant subspaces $\s{M}$ of 
$H_d^2$ 
which contain a polynomial.  In \cite{MT} the theorem is proved for 
certain $\s{M}$ of finite codimension in spaces with complete NP kernels $k$ such 
that $k_\lambda(\lambda) \sra \infty$ as $\lambda \sra \partial B_d$.  

It turns out that vector-valued analogs of Theorem~\ref{Theorem 1.1} and 
Theorem~\ref{Theorem 
1.2} 
are true.  Before we can explain this, we need a few more definitions.  

If $\s{D}$ is a separable complex Hilbert space, then $\s{H}(k, \s{D})$ is the 
space of $\s{D}$-valued $\s{H}(k)$-functions.  It is the set of all analytic 
functions $f:B_d \lra \s{D}$ such that for each $x \in \s{D}$ the function 
$f_x(\lambda) = \left< f(\lambda), x \right>_{\s{D}}$ defines a function in 
$\s{H}(k)$ and such that 

\begin{equation*}
\|f\|^2 = \sum_{n=1}^\infty \, \|f_{e_n}\|^2 < \infty
\end{equation*}

\noindent for some orthonormal basis $\{e_n\}_{n \leq 1}$ of $\s{D}$.  One 
shows that the above expression is independent of the choice of orthonormal 
basis.  In particular, one has for $f \in \s{H}(k), x \in \s{D}$ the function 
$fx: \lambda \lra f(\lambda)x$ is in $\s{H}(k, \s{D})$ and $\|fx\| = \|f\| 
\|x\|_{\s{D}}$.  If $f \in \s{H}(k,\s{D}), x \in \s{D}$, and $\lambda \in B_d$ 
we have $\left< f(\lambda), x \right>_{\s{D}} = \left< f, k_{\lambda}x 
\right>$.  There is an obvious identification of the tensor product $\s{H}(k) 
\otimes \s{D}$ with $\s{H}(k,\s{D})$, where one identifies the elementary 
tensors $f \otimes x$ with the functions $fx$.  Considering the definition of 
the norm in $\s{H}(k, \s{D})$, one may also think of $\s{H}(k, \s{D})$ as a 
direct sum of $\dim \s{D}$ copies of the scalar valued space $\s{H}(k)$.  

Each (scalar valued) multiplier $\varphi \in M(k)$ defines an operator on 
$\s{H}(k, \s{D})$ of the same norm, and we shall also denote this operator by 
$M_\varphi$.  Again, we shall say that a subspace $\s{M}$ of $\s{H}(k, \s{D})$ 
is multiplier invariant if $M_\varphi \s{M} \subset \s{M}$ for each $\varphi 
\in M(k)$.  

Let $\s{D}$ and $\s{E}$ be two separable Hilbert spaces, and let $\phi: B_d 
\lra \s{B}(\s{E}, \s{D})$ be an operator valued analytic function.  For 
$\lambda \in B_d$ and $f \in \s{H}(k, \s{E})$ we define $(\Phi f)(\lambda) = 
\phi(\lambda)f(\lambda)$, then $\Phi f$ is a $\s{D}$-valued analytic function.  
If $\Phi f \in \s{H}(k, \s{D})$ for every $f \in \s{H}(k,\s{D})$, then $\phi$ 
is called an operator-valued multiplier, and one shows that the associated 
multiplication operator $\Phi: \s{H}(k,\s{E}) \lra \s{H}(k,\s{E})$ is bounded.  
It is clear that such multiplication operators $\Phi \in \s{B}(\s{H}(k,\s{E}), 
\s{H}(k,\s{D}))$ intertwine the (scalar) multiplication operators $M_\varphi, 
\varphi \in M(k)$.  It will follow from Lemma~\ref{Lemma 2.2} that for spaces with 
complete NP kernels the multipliers $M(k)$ are dense in $\s{H}(k)$.  It then 
follows from standard arguments that a bounded linear operator $A: 
\s{H}(k,\s{E}) \lra \s{H}(k, \s{D})$ intertwines every $M_\varphi, \varphi \in M(k)$ 
(i.e. $AM_\varphi = M_\varphi A$) if and only if $A=\Phi$ for some 
multiplication 
operator.  

A short calculation shows that for any multiplication operator $\Phi \in 
\s{B}(\s{H}(k, \s{E}), \s{H}(k, \s{D}))$ one has $\Phi^*(k_\lambda x) = k_\lambda 
\phi(\lambda)^* x$ for all $x \in \s{D}, \lambda \in B_d$.  Thus we have 
$\|\phi(\lambda)\|_\s{D} \leq \|\Phi\|$ for all $\lambda \in B_d$ and it 
follows from 
standard arguments that for a.e. $z \in \partial B_d$, $\phi(\lambda)$ converges in 
the strong operator topology to an operator $\phi(z)$ as $\lambda$ approaches $z$ 
nontangentially (for the scalar case see \cite{rudinftotub}, then see 
\cite{rosenblumrovnyak}, pages 81-84 on how to get the operator-valued version).  
Similarly, by applying this reasoning to $\phi(\ol{\lambda})^*$ one sees that also 
$\phi(\lambda)^* \rightarrow \phi(z)^* \sot$ for a.e. 
$z 
\in \partial B_d$ as $\lambda$ approaches $z$ nontangentially.  Actually, the 
limits exist a.e. if the approach is from within certain nonisotropic approach 
regions which for $d>1$ are larger than the standard nontangential approach regions (see 
Section 2 for definitions).  

A multiplication operator $\Phi$ is called {\it inner} if it is a partial isometry as 
an operator $\s{H}(k,\s{E}) \lra \s{H}(k,\s{D})$.  Since partial isometries have 
closed range it is clear that every inner multiplier defines a multiplier invariant 
subspace $\s{M} = \Phi \s{H}(k,\s{E}) \subseteq \s{H}(k,\s{D})$.  Again, it is 
a a remarkable fact that the converse to this theorem is true if $k$ is a 
complete NP kernel.   

\begin{thm} \label{Theorem 1.3}
Let $k$ be a complete NP kernel, let $\s{D}$ be a separable Hilbert space, and let 
$\s{M} \subseteq \s{H}(k, \s{D})$ be a multiplier invariant subspace.  

Then there is an auxiliary Hilbert space $\s{E}$ and an inner multiplication operator 
$\Phi \in \s{B}(\s{H}(k, \s{E}), \s{H}(k, \s{D}))$ such that $\s{M} = \Phi 
\s{H}(k, 
\s{E})$ and $P_\s{M} = \Phi \Phi^*$.  

Furthermore, if $\s{F}$ is another Hilbert space and $\Psi \in \s{B}(\s{H}(k,
\s{F}),$ $\s{H}(k,\s{D}))$ is another inner multiplication operator such that
$\Phi \s{H}(k, \s{E}) = \Psi \s{H}(k, \s{F})$, then there is a partial
isometry $V \in \s{B}(\s{E}, \s{F})$ such that $\phi(\lambda) = \psi(\lambda)
V$ for all $\lambda \in B_d$.  \end{thm}

This theorem is from \cite{MT}, see \cite{arvcurvinv} for the case of 
$H_d ^2$.  Theorem~\ref{Theorem 1.3} implies Theorem~\ref{Theorem 1.1}.  To see this 
we take 
$\s{D} = \mbb{C}$, fix an orthonormal basis $\{e_n\}$ of $\s{E}$ and set 
$\varphi_n(\lambda) = \phi(\lambda) e_n$ for $n \geq 1$ and $\lambda \in B_d$.  With 
this notation it is easy to verify that each $\varphi_n \in M(k) \cap \s{M}$ and 
$P_\s{M} = \Phi \Phi^* = \sum_{n \geq 1} \, M_{\varphi_n} M_{\varphi_n}^* \sot$. 

We note that in the classical Beurling-Lax-Halmos theorem for $H^2 (\mbb{D})$ one may 
take $\s{E} = \s{D}$, but for general complete NP kernels other than the 
Szeg\"{o} kernel that may not be possible if $\dim \s{D} < \infty$.  In fact, one 
may have to take $\s{E}$ to be infinite dimensional even if $\dim \s{D} < \infty$.  
This happens for example for the classical Dirichlet space.  Since the existence of 
these inner multiplication operators is important for our paper, we give a brief 
outline of the proof of Theorem~\ref{Theorem 1.3}.  The proof will explain where the 
space $\s{E}$ come from.  

We already mentioned that the functions $b_n, \; n \geq 1$ in
\eqref{1.1} are multipliers and $P_0 = \sum_{n \geq 1} \, M_{b_n}
M_{b_n}^* \sot$ is the projection onto $\{f \in \s{H}(k) : f(0) = 0 \}$.  Similarly,
if one thinks of the operators $M_{b_n}$ as multipliers on $\s{H}(k,\s{D})$, then
$E_0 = \sum_{n \geq 1} \, M_{b_n} M_{b_n}^* \sot$ is the projection onto $\{f \in
\s{H}(k,\s{D}) : f(0) = 0\}$, and it is easy to see that $Q(A) = \sum_{n \geq 1} \,
M_{b_n} M_{b_n}^*$ defines a completely positive map $\s{B}(\s{H}(k,\s{D})) \lra
\s{B}(\s{H}(k,\s{D}))$.  Now if $\s{M} \subseteq \s{H}(k,\s{D})$ is an multiplier
invariant subspace, then one computes

\begin{equation*}
P_\s{M} - Q(P_\s{M}) = P_\s{M} (I-E_0) P_\s{M} + P_\s{M} Q(I-P_\s{M})P_\s{M} \geq 0.
\end{equation*}

\noindent We set $S = (P_\s{M} - Q(P_\s{M}))^{1/2}, \, \s{E} = (\ker S)^\perp \appa 
\s{H}(k,\s{D})$, and for $\lambda \in B_d, \, x \in \s{D}$, 

\begin{equation*} 
\phi(\lambda)^* x = S(k_\lambda x).
\end{equation*}

\noindent With these definitions one verifies that $\phi$ is an operator valued 
multiplier and that the associated multiplication operator $\Phi$ satisfies $\Phi 
\Phi^* = P_\s{M}$.  

The vector analogue of Theorem~\ref{Theorem 1.2} is that under certain circumstances 
the 
analytic functions associated with inner multiplication operators deserve to be 
called inner functions.  In fact, we shall prove the following theorem.  

\begin{thm} \label{Theorem 1.4}
Let $k$ be a complete NP kernel on $B_d$ and assume that there is a set $\s{P} \appa 
\s{H}(k) \cap C(\ol{B_d})$ which is dense in $\s{H}(k)$ and such that for all $p \in 
\s{P}$ and $z \in \partial B_d, \, \lim_{\lambda \sra z} \, \frac{\|p 
k_\lambda\|}{\|k_\lambda\|} = |p(z)|$.  

Let $\s{E}$ and $\s{D}$ be separable Hilbert spaces, and let $\Phi \in \s{B}(\s{H}(k, 
\s{E}),$ $\s{H}(k, \s{D}))$ be an inner multiplication operator with 
associated 
operator-valued multiplier $\phi: B_d \lra \s{B}(\s{E}, \s{D})$.  

Then for a.e. $z \in B_d$, $\phi(z)$ is a partial isometry with 

\begin{equation}
m = \rank \phi(z) = \sup \{\rank \phi(\lambda): \lambda \in B_d\} = \sup\{\dim 
E_\lambda \s{M} : \lambda \in B_d\}, 
\end{equation}

\noindent where $\s{M} = \Phi \s{H}(k,\s{E})$ and $E_\lambda$ denotes the point 
evaluation map $E_\lambda: \s{H}(k,\s{D}) \lra \s{D}$, $E_\lambda f = f(\lambda)$, $f 
\in \s{H}(k,\s{D})$.  In particular, $m \leq \dim \s{D}$.  
\end{thm}

We shall prove this as Theorem~\ref{Theorem 3.2}.  Of course, as was the case with 
Theorem~\ref{Theorem 1.2}, this theorem applies to the spaces $D_\alpha, \, 0 \leq 
\alpha 
\leq 1$, and $H_d^2$.  In the classical Beurling-Lax-Halmos theorem for 
$H^2(\mbb{D})$ it is also true that the initial space of $\phi(z)$ is a.e. equal to 
one fixed space $\s{K} \appa \s{E} = \s{D}$.  In the general situation of invariant 
subspaces $\s{M}$ of $H_d^2(\s{D}), \, d > 1$, the initial space of $\phi(z)$ may 
vary with $z \in \partial B_d$.  

Section 5 contains our results on the curvature invariant of contractive Hilbert 
modules.  

\vspace{1cm}

\section{Some preliminaries}

Let $k$ be an analytic reproducing kernel on $B_d$ with associated Hilbert space 
$\s{H}(k)$.  We will also assume that $k_0(\lambda) = 1$ for all $\lambda \in B_d$, 
but we will not necessarily assume that $k$ is a complete NP kernel.  

Let $\s{M}$ be multiplier invariant subspace of $\s{H}(k)$.  If $\s{M}$ contains a 
function that does not vanish at $0$, then the function $\varphi_\s{M} = 
\frac{P_\s{M} 1}{\sqrt{P_\s{M} 1(0)}}$ has norm $1$ and solves the extremal problem

\begin{equation} \label{2.1}
\sup \{ \Re f(0) : f \in \s{M}, \, \|f\| = 1\}.
\end{equation}

In fact, if $f \in \s{M}$, then one calculates that $\left< f, \varphi_\s{M} \right> 
= \frac{f(0)}{\varphi_\s{M}(0)}$, and so $\|f\| \geq 
\left| \frac{f(0)}{\varphi_\s{M}(0)}\right|$, which implies the extremal property of 
$\varphi_\s{M}$.  

\begin{lemma} \label{Lemma 2.1}
Let $k$ be a complete NP kernel, let $\s{M} \appa \s{H}(k)$ be a multiplier invariant 
subspace.  Then $P_\s{M} 1 \in M(k)$ and $|(P_\s{M} 1)(0)| \|f\|^2$ for all $f \in 
\s{H}(k)$.  
\end{lemma}

This is proved in \cite{MT}, and it follows immediately from 
Theorem~\ref{Theorem 1.1} or Theorem~\ref{Theorem 1.3}.  In fact, let $\Phi$ be the 
inner 
multiplication operator associated with $\s{M}$ as in \eqref{1.3}, let $\phi$ be the 
associated operator-valued multiplier, take $\s{D} = \mbb{C}$, and set 
$\varphi_n(\lambda) = \left< \phi(\lambda), e_n \right>$, where $\{e_n\}$ is some 
orthonormal basis for the auxiliary space $\s{E}$.  

Then for all $\lambda \in B_d$ we have $(P_\s{M} 1)(\lambda) = \sum_{n \geq 1} \, 
\varphi_n(\lambda) \ol{\varphi_n(0)}$.  It follows that for $f \in \s{H}(k)$, 

\begin{multline*}
\|(P_\s{M} 1)f\|^2 = \|\sum_{n \geq 1} \, \varphi_n \ol{\varphi_n (0)} f \|^2 \leq 
\|\Phi\|^2 \sum_{n \geq 1} \, \|\ol{\varphi_n (0)} f \|^2 \\
\leq \sum_{n \geq 1} \, |\varphi_n (0) |^2 \|f\|^2 = |P_\s{M} 1)(0)| \|f\|^2.  
\end{multline*}

For $\lambda \in B_d$ we define $\s{M}_\lambda = \{f \in \s{H}(k): f(\lambda) = 0\}$.  
Then each $\s{M}$ is a multiplier invariant subspace of $\s{H}(k)$ with 
$\s{M}_\lambda = \{k_\lambda\}^\perp$.  Thus, for $\lambda \neq 0$, one obtains 

\begin{equation} \label{2.2}
\varphi_\lambda(z) = \varphi_{\s{M}_\lambda} (z) = \frac{1 - k_\lambda 
(z)/k_\lambda(\lambda)}{\sqrt{1-1/k_\lambda(\lambda)}}.
\end{equation}

We shall refer to $\varphi_\lambda$ as the {\it one point extremal function}.  Note 
that \eqref{2.1} implies that if $k$ is a complete NP kernel, then all one point 
extremal functions $\varphi_\lambda, \, \lambda \neq 0$ contractive multipliers.  

\begin{lemma} \label{Lemma 2.2}
Let $k$ be such that for each $\lambda \in B_d \backslash \{0\}$ the one point 
extremal function $\varphi_\lambda$ is a contractive multiplier on $\s{H}(k)$, let 
$\s{D}$ be a separable Hilbert space, and let $f \in \s{H}(k, \s{D})$.  

Then 

\begin{enumerate}
%here!
\item \label{a}
for each $\lambda \in B_d$, $k_\lambda \in M(k)$ and $\|k_\lambda 
\|_M \leq 2 
k_\lambda(\lambda)$,  
\item \label{b}
for each $\lambda \in B_d$, $\|f(\lambda)\|^2_\s{D} \leq 
\displaystyle{\frac{\|k_\lambda 
f\|^2}{\|k_\lambda\|^2}} \leq 2 \Re \left< f, k_\lambda f \right> - \|f\|^2$, 
\item \label{c}
the function $F: B_d \lra \mbb{C}, \, F(\lambda) = \left< f, k_\lambda 
f 
\right>$ is analytic on $B_d$.
\end{enumerate}
\end{lemma}

\begin{proof}
For $\lambda = 0$, \eqref{a} and \eqref{b} are clear since $k_0 \equiv 1$.  For 
$\lambda \neq 0$, \eqref{a} follows from \eqref{2.2} since $\varphi_\lambda \in M(k)$ 
and $\|\varphi_\lambda\| \leq 1$.  Furthermore, the hypothesis also implies that 
$\|f\|^2 - \|\varphi_\lambda f\|^2 \geq 0$.  After a short calculation this leads to 
the right inequality of \eqref{b}.  To see the left inequality in \eqref{b}, note 
that $k_\lambda f(\lambda) : z \mapsto k_\lambda (z) f(\lambda)$ defines a function 
in $\s{H}(k,\s{D})$ which is orthogonal to $(k_\lambda f) - k_\lambda f(\lambda)$.  
Hence $\|k_\lambda f\|^2 = \|k_\lambda f - k_\lambda f(\lambda) \|^2 \geq \|k_\lambda 
\|^2 \|f(\lambda)\|_\s{D}^2$.  

We now prove \eqref{c}.  If $f \in \s{H}(k,\s{D})$, $\varphi \in M(k)$, and $x \in 
\s{D}$, then $\left< f, k_\lambda \varphi x \right> = \left< (M_\varphi^* 
f)(\lambda), x \right>_\s{D}$ is an analytic function in $\lambda \in B_d$.  Hence if 
$\s{L} \appa \s{H}(k,\s{D})$ is the set of finite linear combinations of elements of 
the form $\varphi x$, $\varphi \in M(k)$, $x \in \s{D}$, then for each $f \in 
\s{H}(k,\s{D})$ and $g \in \s{L}$ the function $\lambda \mapsto \left< f, k_\lambda g 
\right>$ is analytic in $B_d$.  Finite linear combinations of the functions 
$k_\lambda$ are dense in $\s{H}(k)$, hence it follows from \eqref{a} that $M(k)$ is 
dense in $\s{H}(k)$, and so $\s{L}$ is dense in $\s{H}(k,\s{D})$.  The unit ball in 
$\s{H}(k)$ is a normal family, thus the uniform boundedness principle implies that 
for each compact subset $K \appa B_d$ there is $C_K$ such that $k_\lambda(\lambda) = 
\|k_\lambda\|^2 \leq C_K$.  This implies that for each compact set $K \appa B_d$, 
$\lambda \in K$, $f \in \s{H}(k, \s{D})$, and $g \in \s{L}$, we have by \eqref{a}

\begin{equation*}
\left| \left< f, k_\lambda f \right> - \left< f, k_\lambda g \right> \right| \leq 
\|f\| \|k_\lambda \|_M \|f-g\| \leq 2 C_K \|f\| \|f-g\|,
\end{equation*}

\noindent i.e. for each $f \in \s{H}(k,\s{D}), \, F(\lambda) = \left< f, k_\lambda 
f \right>$ is analytic as it is a local uniform limit of analytic functions.  This 
concludes the proof of the lemma.  
\end{proof}

\noindent This lemma has a number of repercussions for the regularity of the 
functions in $\s{H}(k)$.  

\begin{prop} \label{Proposition 2.3}
Let $k$ be such that for each $\lambda \in B_d \backslash \{0\}$ the one point 
extremal function $\varphi_\lambda$ is a contractive multiplier on $\s{H}(k)$.  
Then $\s{H}(k)$ is contractively contained in $H^2(\partial B_d)$.  In fact, for 
every $z \in \partial B_d$, the slice function $f_z, \, f_z(\zeta) = f(\zeta z), \, 
\zeta \in \mbb{D}$, is in $H^2(\mbb{D})$, and satisfies $\|f_z\|_{H^2} \leq 
\|f\|$.  
\end{prop}

\begin{proof} We use the scalar version of Lemma~\ref{Lemma 2.2}~\eqref{1}.  Hence
for each $f \in \s{H}(k)$ and $\lambda \in B_d$ we have $|f(\lambda)|^2 \leq
u(\lambda) = 2 \Re \left< f, k_\lambda f \right> - \|f\|^2$.  As before, let $\sigma$ 
denote the rotationally invariant probability measure on $\prtl B_d$.  We fix $0 
< r < 1$ and integrate over $\partial B_d$ and obtain 

\begin{equation*}
\int_{\partial B_d} |f(rz)|^2 d\sigma(z) \leq \int_{\prtl B_d} u(rz) d\sigma(z) = 
u(0) = \|f\|^2,
\end{equation*}

\noindent since the integrand on the right is the real part of an analytic function.  
We now take the supremum over $0 < r < 1$ and obtain $\|f\|_{H^2(\prtl B_d)} \leq 
\|f\|$ for all $f \in \s{H}(k)$.  

Furthermore, if $z \in \prtl B_d$, then $u_z(\zeta) = u(\zeta z), \, \zeta \in 
\mbb{D}$, defines a positive harmonic function in the unit disc $\mbb{D} \appa 
\mbb{C}$.  Thus $|f_z(\zeta)|^2 = |f(\zeta z)|^2 \leq u_z(\zeta)$, hence 
$\|f_z\|^2_{H^2} \leq u_z(0) = \|f\|^2$.  
\end{proof}

Functions in $H^2(\prtl B_d)$ have a.e. limits from within certain approach 
regions 
that contain the standard nontangential approach regions (see 
\cite{rosenblumrovnyak}).  For $\alpha > 1$ and $z \in \prtl B_d$, define 
$\Omega_\alpha (z)$ to be the set of all $\lambda \in B_d$ such that $|1 - \left< 
\lambda, z \right>_d| < \displaystyle{\frac{\alpha}{2}} (1-|\lambda|^2)$.  We say 
that a function $f : B_d \lra \mbb{C}$ has a $K$-limit $A$ at $z \in \prtl B_d$, 
$(K-\lim f)(z) = A$, if for every $\alpha > 1$ and for every sequence $\{\lambda_n\} 
\appa \Omega_\alpha (z)$ that converges to $z$, we have $f(\lambda_n) \sra A$ as $n 
\sra \infty$.  

Let $k$ be a reproducing kernel as in Proposition~\ref{Proposition 2.3}.  It 
is 
well-known that every function in $H^2(\prtl B_d)$ has finite $K$-limits at a.e. 
every point $z \in \prtl B_d$ (\cite{rudinftotub}), hence the same is true for every 
$f \in \s{H}(k)$.  Furthermore, if $\s{D}$ is a separable Hilbert space, and $f \in 
\s{H}(k,\s{D})$, then the arguments given in \cite{rosenblumrovnyak} on page 84 show 
that for a.e. $z \in \prtl B_d$ there is an $f(z) \in \s{D}$ such that $f(z)$ is the 
$K$-norm-limit of $f(\lambda)$ at $z$.  

\begin{prop} \label{Proposition 2.4}
Let $k$ be such that for each $\lambda \in B_d \backslash \{0\}$ the one point 
extremal 
function $\varphi_\lambda$ is a contractive multiplier on $\s{H}(k)$, and assume that 
there is a set $\s{P} \appa \s{H}(k) \cap C(\ol{B_d})$ which is dense in $\s{H}(k)$ 
and such that for all $p \in \s{P}$ and $z \in \prtl B_d, \, \lim_{\lambda \sra z} \, 
\ds{\frac{\|p k_\lambda\|}{\|k_\lambda\|}} = |p(z)|$.  Let $\s{D}$ be a separable 
Hilbert space.  Then for every $f \in \s{H}(k, \s{D})$ we have 

\begin{equation*} 
K-\lim \frac{\|f k_\lambda\|}{\|k_\lambda\|} = \|f(z)\|_\s{D} \quad \text{for} \; 
a.e. z \in \prtl B_d.
\end{equation*}

\end{prop}

\begin{proof}
Because of the hypothesis and Lemma~\ref{Lemma 2.2}, we can use standard techniques.  
We briefly outline the details of the proof.  

If $f \in \s{H}(k,\s{D}), \, \alpha > 1$, we define the maximal function

\begin{equation*}
M_\alpha f(z) = \sup \{\frac{\|f k_\lambda\|}{\|k_\lambda\|}: \lambda \in 
\Omega_\alpha (z)\}.
\end{equation*}

The right hand side in Lemma~\ref{Lemma 2.2}~\eqref{b} is positive and the real part 
of an analytic function (i.e. it is pluriharmonic), hence it can be represented as 
the invariant Poisson integral of a positive measure $\mu$ on $\prtl B_d, \, 
P\mu(\lambda) = 2 \Re \left< f, k_\lambda f\right> - \|f\|^2$ (see 
\cite{rudinftotub}).  Furthermore, we note that $\|\mu\| = P\mu(0) = \|f\|^2$, and 
that for all $\alpha > 1$, the $\Omega_\alpha$-maximal function of $P\mu$ satisfies a 
weak-type estimate with constant $C_\alpha$ (see \cite{rudinftotub}).  Hence by 
Lemma~\ref{Lemma 2.2}~\eqref{b} we obtain for all $\alpha > 1, \, \epsilon > 0$, and 
$f \in \s{H}(k,\s{D})$ the weak-type estimate

\begin{equation*}
\sigma(\{z \in \prtl B_d : M_\alpha f(z) > \epsilon \}) \leq C_\alpha 
\frac{\|f\|^2}{\epsilon^2}.
\end{equation*}

Next, let $\s{P}' \appa \s{H}(k,\s{D})$ be the set of all finite linear combinations 
of the form $px$, where $p \in \s{P} \appa \s{H}(k)$ and $x \in \s{D}$.  Then 
$\s{P}'$ is dense in $\s{H}(k,\s{D})$.  We shall first show that $\lim_{\lambda \sra 
z} \, \ds{\frac{\|(p-p(\lambda))k_\lambda\|}{\|k_\lambda\|}} = 0$ for all $p \in 
\s{P}'$ and $z \in \prtl B_d$.  

Let $z \in \prtl B_d$, and note that if $q \in \s{P}$, then 

\begin{equation*} 
\frac{\|(q-q(\lambda))k_\lambda \|^2}{\|k_\lambda\|^2} = 
\frac{\|qk_\lambda\|^2}{\|k_\lambda\|^2} - |q(\lambda)|^2 \sra 0 \; \text{as} 
\; 
\lambda 
\sra z, 
\end{equation*}

\noindent because $q \in \s{P}$ and $q$ is continuous at $z$.  Now let $p = 
\sum_{i=1}^n \, p_i x_i$, where $p_i \in \s{P}$ and $x_i \in \s{D}$, then 

\begin{equation*}
\frac{\|(p-p(\lambda))k_\lambda\|}{\|k_\lambda\|} \leq \sum_{i=1}^n \, \|x_i\|_\s{D} 
\frac{\|(p_i - p_i(\lambda))k_\lambda\|}{\|k_\lambda\|} \sra 0 \; \text{as} \; 
\lambda 
\sra z.
\end{equation*}

Finally, let $f \in \s{H}(k,\s{D})$.  Then $f$ has $K$-limit $f(z)$ at a.e. $z \in 
\prtl B_d$ and $\ds{\frac{\|(f-f(\lambda))k_\lambda\|^2}{\|k_\lambda\|^2}}$ equals 
$0$ 
at $\sigma$ a.e. $z \in \prtl B_d$.  Using Lemma~\ref{Lemma 2.2}~\eqref{b}, we see 
that for every $p \in \s{P}'$ we have for all $\lambda \in B_d$, 

\begin{multline*} 
\frac{\|(f-f(\lambda))k_\lambda\|}{\|k_\lambda\|} \leq 
\frac{\|(f-p)k_\lambda\|}{\|k_\lambda\|} + 
\frac{\|(p-p(\lambda))k_\lambda\|}{\|k_\lambda\|} + \|p(\lambda) - 
f(\lambda)\|_\s{D} \\
\leq 2 \frac{\|(f-p)k_\lambda\|}{\|k_\lambda\|} + 
\frac{\|(p-p(\lambda))k_\lambda\|}{\|k_\lambda\|}.
\end{multline*}

\noindent Hence for $z \in \prtl B_d$, we obtain for every $\alpha > 1$,

\begin{equation*}
\limsup_{\ds{\stackrel{\lambda \sra z}
{\lambda \in \Omega_\alpha}}} \frac{\|(f-f(\lambda))k_\lambda\|}{\|k_\lambda\|} \leq 
2 M_\alpha (f-p)(z),
\end{equation*}

\noindent and so the weak-type estimate implies that for every $\epsilon > 0$, we 
have 
for every $p \in \s{P}'$,

\begin{equation*}
\sigma(\{z \in \prtl B_d: \limsup_{\substack{\lambda \sra z \\
\lambda \in 
\Omega_\alpha}} > \epsilon\}) \leq 4 C_\alpha \frac{\|f-p\|^2}{\epsilon^2}.
\end{equation*}

\noindent Since $\s{P}'$ is dense in $\s{H}(k,\s{D})$ the result follows.  
\end{proof}

\vspace{1cm}

\section{Inner multiplication operators and inner multipliers.}

\noindent As in Section 2, in this section $k$ will denote an analytic 
reproducing 
kernel on $B_d$ with $k_0 \equiv 1$.  

\begin{lemma} \label{Lemma 3.1}
Let $\s{D}$, $\s{E}$ be separable Hilbert spaces, and let $\Phi$ $\in 
\s{B}(\s{H}(k,\s{E}),$ $\s{H}(k,\s{D}))$ be a multiplication operator with 
associated operator-valued multiplier $\phi, \, \phi(\lambda) \in \s{B}(\s{E}, 
\s{D}), \, \lambda \in B_d$.  For $\lambda \in \ol{B_d}$ let $\rank 
\phi(\lambda) = \dim \ran \phi(\lambda)$, and set 

\begin{equation*}
m = \sup \{\rank \phi(\lambda) : \lambda \in B_d\}.
\end{equation*}

Then $\rank \phi(\lambda) = m$ on $B_d \backslash E$, where $E$ is at most a 
countable union of zero varieties of nonzero bounded analytic functions in 
$B_d$ and $\rank \phi(z) = m$ for $\sigma$ a.e. $z \in \prtl B_d$.  
\end{lemma}

\begin{proof}
First note that if $T_n, \, T \in \s{B}(\s{E}, \s{D})$ such that $T_n \sra T 
\sot$, then $\rank T \leq \liminf \rank T_n$.  In fact, for the proof we may 
assume that $\rank T_n = r < \infty$ for all $n$.  Suppose $\rank T > r$.  
Then let $\{Tf_j\}_{j=1}^{r+1}$ be an orthonormal set in the range of $T, \, 
f_j \in \s{D}$.  Then $d_n = \det (\left< T_n f_j, T f_k \right>_\s{D}) = 0$ 
for each $n$ since $\rank T_n = r < r+1$.  But this leads to a contradiction 
since $d_n \sra \det (\left< T f_j, T f_k \right> _\s{D}) = 1$.  Thus, at each 
point $z \in \prtl B_d$ where the $K$-limit of $\phi(\lambda)$ exists in the 
strong operator topology, we have $\rank \phi(z) \leq m$.  

Now assume that $1 \leq m < \infty$.  Then there is a $\lambda_0 \in B_d$ such 
that $\rank \phi(\lambda_0) = m$.  Let $\{e_n\}_{n \geq 1}^m$ be an 
orthonormal basis for $\ker \phi(\lambda_0)^\perp \appa \s{E}$, and let 
$\{d_k\}_{k \geq 1}^m$ be an orthonormal basis for $\ran \phi(\lambda_0) \appa 
\s{D}$.  

We define $D(\lambda) = \det \left[ (\left< \phi(\lambda) e_n, 
d_k\right>_\s{D})_{1 \leq n,k \leq m} \right]$.  Then $D$ is a bounded 
analytic function in $B_d$ with $D(\lambda_0) \neq 0$.  It is clear that $m 
\leq \rank \phi(\lambda)$ whenever $D(\lambda) \neq 0, \, \lambda \in B_d$, 
but since $m$ was the supremum of $\rank \phi(\lambda)$ for $\lambda \in B_d$ 
we actually get $m = \rank \phi(\lambda)$ whenever $D(\lambda) \neq 0, \, 
\lambda \in B_d$.  

Furthermore, since the determinant is a polynomial in its entries it is clear 
that the $K$-limit of $D(\lambda)$ exists, is nonzero at a.e. $z \in \prtl 
B_d$, and equals $D(z) = \det \left[ (\left< \phi(z) e_n, d_k 
\right>_\s{D})_{1 \leq n,k \leq m}\right]$.  Hence $\rank \phi(z) \geq m$ for 
$\sigma$ a.e. $z \in \prtl B_d$.  We already explained the other inequality, 
thus this proves the lemma when $m < \infty$.  

If $m = \infty$, then for any integer $s>0$ we can find $\lambda_s$ such that 
$\rank \phi(\lambda_s) \geq s$.  Thus, as above, we obtain a bounded analytic 
function $D_s(\lambda)$ with $D_s (\lambda_s) \neq 0$.  The boundary value 
function of $D_s$ is not identically zero, hence $\rank \phi(z) \geq s$ for 
a.e. $z \in \prtl B_d$.  This implies that $\rank \phi(z) = \infty$ for a.e. 
$z \in \prtl B_d$.  It also follows that $\rank \phi(\lambda) = \infty$ for 
all $\lambda \in B_d, \, \lambda \notin \bigcap_{k=1}^\infty 
\bigcup_{s=k}^\infty \, Z(D_s)$.  
\end{proof}

For $\lambda \in B_d$, let $E_\lambda : \s{H}(k,\s{D}) \lra \s{D}, \, 
E_\lambda f = f(\lambda)$, let $\Phi \in \s{B}(\s{H}(k,\s{E})$ 
$\s{H}(k,\s{D}))$ be an inner multiplication operator with associated 
operator-valued multiplier $\phi$, and let $\s{M} = \Phi \s{H}(k,\s{E})$.  
Then for all $\lambda \in B_d$, $\ran \phi(\lambda) = \{\phi(\lambda) y : y 
\in \s{E} \} = \{\phi(\lambda) f(\lambda) : f \in \s{H}(k,\s{E}) \} = 
\{E_\lambda(\Phi f) : f \in \s{H}(k,\s{E}) \} = E_\lambda \s{M}$.  

\begin{thm} \label{Theorem 3.2}
Let $k$ be such that for each $\lambda \in B_d \backslash \{0\}$ the one point 
extremal function $\varphi_\lambda$ is a contractive multiplier on $\s{H}(k)$ 
and assume that there is a set $\s{P} \appa \s{H}(k) \cap C(\ol{B_d})$ which 
is dense in $\s{H}(k)$ and such that for all $p \in \s{P}$ and $z \in B_d$, 
$\lim_{\lambda \sra z} \ds{\frac{\|p k_\lambda\|}{\|k_\lambda\|}} = |p(z)|$.  

Let $\s{E}$ and $\s{D}$ be separable Hilbert spaces, and let $\Phi \in 
\s{B}(\s{H}(k,\s{E}),$ $\s{H}(k,\s{D}))$ be an inner multiplication operator 
with 
associated operator-valued mulitplier $\phi$, and let $m$ $ = \sup \{\rank 
\phi(\lambda) : \lambda \in B_d\}$ $ = \sup $ $\{\dim E_\lambda \s{M} : 
\lambda 
\in B_d 
\}$.  

Then for $\sigma$ a.e. $z \in \prtl B_d$, $\phi(z)$ is a partial isometry with 
$\rank \phi(z) = m$.  
\end{thm}

\begin{proof}
The statment about the rank follows from Lemma~\ref{Lemma 3.1}.  Let $z \in 
\prtl B_d$ be such that the $K$-limit of $\phi(\lambda)$ exists at $z$ in the 
strong operator topology.  We have to show that $\phi(z)^*$ is an isometry on 
$\ran \phi(z)$.  Since $\|\phi(z)^*\| \leq 1$ it suffices to show that 
$\|\phi(z)^* \phi(z) y\|_\s{D} \geq \|\phi(z)y\|_\s{D}$ for all $y \in \s{E}$ 
with $\phi(z)y \neq 0$.  

Let $\s{M} = \ran \Phi \appa \s{H}(k,\s{D})$, then $P_\s{M} = \Phi \Phi^*$ and 
$\s{M}$ is a multiplier invariant subspace.  Note that for $\lambda \in B_d$ 
and $x \in \s{D}$ we have 

\begin{equation} 
\|P_\s{M} k_\lambda x\| = \sup \{|\left< f(\lambda), x \right>_\s{D} | : f \in 
\s{M}, \, \|f\| \leq 1 \},
\end{equation}

\noindent because for $f \in \s{M}$, $|\left< f(\lambda), x\right>_\s{D}| = 
|\left< f, k_\lambda x \right>| = |\left< f, P_\s{M} k_\lambda x \right> | 
\leq \|f\| \|P_\s{M} k_\lambda x \|$ with equality if $f = P_\s{M} k_\lambda 
x$.  

Hence if $f \in \s{M}$ is nonzero, if $\lambda \in B_d$, then $k_\lambda f \in 
\s{M}$ and for $x \in \s{D}$

\begin{multline*}
\|\phi(\lambda)^* x\|_\s{D} = \frac{\|k_\lambda \phi(\lambda)^* x 
\|}{\|k_\lambda\|} = \frac{\|\Phi^* (k_\lambda x) \|}{\|k_\lambda\|} = 
\frac{P_\s{M} (k_\lambda x) \|}{\|k_\lambda\|} \\
\geq \frac{|\left< (k_\lambda f)(\lambda), x \right>_\s{D} |}{\|k_\lambda\| \,  
\|k_\lambda f \|} = \frac{|\left< f(\lambda), x \right>_\s{D}|}{\|k_\lambda 
f\|/\|k_\lambda\|}.
\end{multline*}

Now let $y \in \s{E}$ with $\phi(z)y \neq 0$.  Then $f = \Phi y$, $f(\lambda) 
= \phi(\lambda) y$, is a nonzero function in $\s{M}$ with $K$-limit $f(z) = 
\phi(z)y$ as $\lambda \sra z$.  Hence Proposition~\ref{Proposition 2.4} and 
the above imply $\|\phi(z)^* x\|_\s{D} \geq \ds{\frac{|\left< \phi(z)y, 
x\right>_\s{D}|}{\|\phi(z)y\|}}$.  This concludes the proof since we may take 
$x =  \phi(z)y$.  
\end{proof}

\vspace{1cm}

\section{$\s{U}$-invariant complete NP kernels.}

In this section we shall verify that the hypothesis of 
Theorem~\ref{Theorem 3.2} is 
satisfied for many complete NP kernels that are invariant under unitary maps.  
If $k$ is an analytic reproducing kernel on $B_d$ that is invariant under 
every unitary map $U: \mbb{C}^d \lra \mbb{C}^d$, i.e. $k_{U\lambda}(Uz) = 
k_\lambda(z)$ for all $\lambda, z \in B_d$, then one can show that $k_\lambda 
(z) = f(\left< \lambda,z\right>_d)$ for some function $k$ of type $f(x) = 
\sum_{n=0}^\infty \, a_n x^n, \, a_n \geq 0$.  Such a kernel will be a 
complete NP kernel if and only if there exists a sequence $\{b_n\}_n, \, b_n 
\geq 0$ such that $f(x) = \ds{\frac{1}{1-\sum_{n=1}^\infty \, b_n x^n}}$.  
Note that 
$a_1 = b_1$.  Finite linear combinations of the kernels $k_\lambda, \, \lambda 
\in B_d$, are dense in $\s{H}(k)$ and the evaluations of partial derivatives 
of $0$ are continuous linear functionals, hence if $a_1 = 0$, it would follow 
that the coordinate functions $z_i(\lambda) = \lambda_i, \, i = 1, \ldots, d$ 
are not in $\s{H}(k)$.  On the other hand, if $a_1 = b_1 \neq 0$, then each 
$z_i \in M(k)$.  This follows because we have $1-\ds{\frac{1}{k_\lambda(w)}} = 
b_1 \sum_{i=1}^d \, w_i \ol{\lambda}_i +$ higher order terms.  Thus 
$\sum_{i=1}^d \, M_{z_i} M_{z_i}^* \leq (1/b_1)I$ because, as we already 
mentioned, each of the functions $b_n$ in the representation \eqref{1.1} is a 
multiplier with $\sum_{n \geq 1} \, M_{b_n} M_{b_n}^* = P_0 \leq I$.  It 
follows from this, or it is easy to see anyway, that the hypothesis $b_1 \neq 
0$ implies that $a_n > 0$ for all $n$.  

Thus in this section we will assume that 

\begin{equation} \label{4.1}
k_\lambda(z) = \sum_{n=0}^\infty \, a_n (\left< z,\lambda \right>_d)^n = 
\frac{1}{1-\sum_{n=1}^\infty \, b_n (\left< z, \lambda \right>_d)^n}
\end{equation}

\noindent where $a_n, \, b_n \geq 0$, $a_0 = 1$, and $a_1 = b_1 > 0$.  In 
particular, $k$ is a complete NP kernel, and the space $\s{H}(k)$ contains the 
polynomials.  

At this point we should mention that if only the sequence $\{a_n\}$ is given, 
then it may be difficult to determine whether $k$ is a complete NP kernel, 
i.e. for which $\{a_n\}$ it follows that a sequence $\{b_n\}$ can be found 
such that $b_n \geq 0$ for each $n$ and such that \eqref{4.1} holds.  However, 
it was pointed out in \cite{shapiroshields1} that if $a_{n+1}/a_n$ increases 
to $1$, then the existence of nonnegative $\{b_n\}$ follows by a theorem of 
Hardy, \cite{hardy}.  On the other hand, if the sequence $\{b_n\}$ is given, 
$b_n \geq 0$, one always obtains a complete NP kernel.  

In order to compute the norm of polynomials we need to recall multiindex 
notation.  Let $k = (k_1, k_2, \ldots, k_d)$ be a multiindex of nonnegative 
integers, then $|k| = k_1 + k_2 + \cdots + k_d$, $k! = k_1! k_2! \cdots k_d!$, 
and for $\lambda = (\lambda_1, \lambda_2, \ldots, \lambda_d) \in \mbb{C}^d$, 
$\lambda^k = \lambda_1^k \lambda_2^k \cdots \lambda_d^k$, and the multinomial 
formula implies that for $z, \lambda \in B_d$ and $n \geq 0$

\begin{equation*}
\left< z,\lambda\right>_d^n = \sum_{|k| = n} \, \frac{|k|!}{k!} z^k 
\ol{\lambda}^k.
\end{equation*}

Thus $k_\lambda(z) = \sum_k \, a_{|k|} \ds{\frac{|k|!}{k!}} z^k \ol{w}^k$, 
where 
the 
sum is taken over all multiindices $k$ with entries in the nonnegative 
integers.  Since $k_\lambda(z) = \left< k_\lambda, k_z \right>$ it follows 
that monomials in $\s{H}(k)$ are mutually orthogonal and 

\begin{equation*}
\|z^k\|^2 = \frac{k!}{a_{|k|} |k|!}.
\end{equation*}

If $f_n(z) = \sum_{|k| = n} \, c_k z^k$ is a homogeneous polynomial of degree 
$n$, then $\|f_n\|^2 = \ds{\frac{1}{a_n}} \sum_{|k| = n} \, 
\ds{\frac{k!}{|k|!}} |c_k|^2$, and it follows that an analytic function $f$ on 
$B_d$ with homogeneous expansion $f(z) = \sum_{n=0}^\infty \, f_n(z)$ is in 
$\s{H}(k)$ if and only if $\|f\|^2 = \sum_{n=0}^\infty \, \|f_n\|^2 < \infty$, 
and the polynomials are dense in $\s{H}(k)$.  

We need a few technical lemmas.  

\begin{lemma} \label{Lemma 4.1}
Let $\s{D}$ be a separable Hilbert space, and for $n \in \mbb{N}$ let 
$F_n(e^{it}) =\ds{\frac{1}{n}} \ds{\frac{\sin^2 ((n+1)t/2)}{\sin^2(t/2)}}$ be 
the 
Fejer kernel.  

If $\varphi \in M(k), \, n \in \mbb{N}$, and $p_n(z) = \int_0^{2\pi} \, 
\varphi(e^{it} z) F_n (e^{it}) \, dt$, then $M_{p_n} \sra M_\varphi \wot$ in 
$\s{H}(k)$ or $\s{H}(k, \s{D})$.  

In particular, if a subspace $\s{M}$ is $M_{z_i}$ invariant for each $i = 1, 
\ldots, d$, then it is multiplier invariant.  
\end{lemma}

\begin{proof}
For $t \in \mbb{R}$, $f \in \s{H}(k, \s{D})$ let $f_t(z) = f(e^{it} z)$.  One 
checks that $f_t \in \s{H}(k, \s{D})$, $\|f_t\| = \|f\|$, and $f_t \sra 
f_{t_0}$ in norm as $t \sra t_0$.  Hence if $\varphi \in M(k)$, then $\varphi 
f = (\varphi f_{-t})_t$, hence $\varphi_t \in M(k)$ with $\|\varphi_t\|_M = 
\|\varphi\|_M$, and $\varphi_t f \sra \varphi_{t_0} f$ in norm as $t \sra 
t_0$.  Hence for each $n \in \mbb{N}$ the integral $p_n f = \int_0^{2\pi} \, 
\varphi_t f F_n(e^{it}) \, dt$ converges in the norm of $\s{H}(k, \s{D})$, and 
we have $\|p_n f\| \leq \|\varphi\|_M \|f\|$.  

The lemma follows from this and the fact that $p_n (\lambda) \sra 
\varphi(\lambda)$ for each $\lambda \in B_d$.  
\end{proof}

\begin{lemma} \label{Lemma 4.2}
Suppose $k$ satisfies \eqref{4.1}.  Let $p$ be a homogeneous polynomial of 
degree $n$.  Then

\begin{equation}
\sum_{i=1}^d \, \|z_i p\|^2 - \|p\|^2 = \left( \frac{a_n}{a_{n+1}} 
\frac{n+d}{n+1} -1 \right) \|p\|^2.
\end{equation}

\noindent Hence, if $a_n/a_{n+1} \sra 1$ as $n \sra \infty$, then 
$\sum_{i=1}^d \, M_{z_i}^*M_{z_i} - I$ is a compact operator on $\s{H}(k)$.  
\end{lemma}

\begin{proof}
For $1 \leq i \leq d$ and any multiindex $k = (k_1, \ldots, k_d)$, we obtain 
$\|z_i z^k\|^2 = \ds{\frac{a_{|k|}}{a_{|k| +1}} \frac{k_i + 1}{|k| + 1}} 
\|z^k\|^2$.  Thus, if $p$ is a homogeneous polynomial of degree $n \geq 0$, 
then $\sum_{i=1}^d \, \|z_i p\|^2 = \ds{\frac{a_n}{a_{n+1}} \frac{n+d}{n+1}} 
\|p\|^2$.  
\end{proof}

\begin{thm} \label{Theorem 4.3}
Let $\s{E}$ and $\s{D}$ be separable Hilbert spaces, and let $k$ be a complete 
NP kernel that satisfies the hypothesis \eqref{4.1}, $k_\lambda (\lambda) 
\sra \infty$ as $|\lambda| \sra 1$, and $\lim_{n \sra \infty} \, a_n/a_{n+1} 
= 
1$. 

If $\Phi \in \s{B}(\s{H}(k, \s{E}), \s{H}(k, \s{D}))$ is an inner 
multiplication operator with associated operator-valued multiplier $\phi$, and 
$m = \sup \{\rank \phi(\lambda) : \lambda \in B_d \}$, then for $\sigma$ a.e. 
$z \in B_d$, $\phi(z)$ is a partial isometry with $\rank \phi(z) = m$.  

Furthermore, Theorem~\ref{Theorem 1.3} applies to every subspace $\s{M} \appa 
\s{H}(k, \s{D})$ that is invariant for every $M_{z_i}, \, i = 1, \ldots, d$.  
In this case, $m = \sup \{ \dim$ $E_\lambda \s{M} : \lambda \in B_d \}$.  
\end{thm}

\begin{proof}
The statement of the last paragraph follows from Lemma~\ref{Lemma 4.1} and 
Theorem~\ref{Theorem 3.2}.  We shall show that $\lim_{\lambda \sra w} \, 
\ds{\frac{\|p k_\lambda\|}{\|k_\lambda\|}} = |p(w)|$ for every $w \in \prtl 
B_d$ and every polynomial $p$.  

If $p$ is a polynomial and $w \in \prtl B_d$, then there are polynomials $q_i, 
\, i = 1, \ldots, d$ such that $p(z) - p(w) = \sum_{i=1}^d \, (z_i - w_i) 
q_i(z)$.  Then for $\lambda \in B_d$, 

\begin{multline*}
\frac{\|(p-p(w))k_\lambda\|}{\|k_\lambda\|} \leq \sum_{i=1}^d \, \|q_i\|_M 
\frac{\|(z_i - w_i)k_\lambda \|}{\|k_\lambda\|} \leq C\left[ \sum_{i=1}^d \, 
\frac{\|(z_i - w_i)k_\lambda\|^2}{\|k_\lambda\|^2} \right]^{1/2} \\
= C \left[ \sum_{i=1}^d \, \frac{z_i k_\lambda \|^2}{\|k_\lambda\|^2} - 2 \Re 
\lambda_i \ol{w}_i + |w_i|^2 \right]^{1/2} \\
= C \left[ \sum_{i=1}^d \, \|z_i \frac{k_\lambda}{\|k_\lambda\|}\|^2 - 2 \Re 
\left< \lambda, w \right> + 1 \right]^{1/2}.
\end{multline*}

Now the hypothesis $k_\lambda(\lambda) \sra \infty$ and the density of the 
polynomials implies that $k_\lambda/\|k_\lambda\| \sra 0$ weakly as 
$|\lambda| \sra 1$, and by Lemma~\ref{Lemma 4.2} we have that $\sum_{i=1}^d \, 
M_{z_i}^* M_{z_i} - I$ is compact, so $\sum_{i=1}^d \, \|z_i 
(k_\lambda/\|k_\lambda\|)\|^2 \sra 1$ as $|\lambda| \sra 1$.  Hence 
$\ds{\frac{\|(p-p(w)) k_{\lambda} \|}{\|k_\lambda\|}} \sra 0$ as $\lambda \sra 
w$.  
\end{proof}

%font?
\begin{exam} \label{Example 4.4} 
Let $a_n = (n+1)^{-\alpha}, \, \alpha \geq 0$.  Then the corresponding kernel 
$k$ is a complete NP kernel.  This follows from the theorem of Hardy that we 
already mentioned after \eqref{4.1}.  For $d=1$ the spaces $\s{H}(k) = 
D_\alpha$ are weighted Dirichlet spaces with $D= D_1$ being the classical 
Dirichlet space.  We note that for $0 \leq \alpha \leq 1$ the coefficients 
$a_n$ satisfy the hypothesis of the theorem.  However, if $\alpha > 1$, then 
$k_\lambda(\lambda)$ stays bounded as $|\lambda| \sra 1$.  

We shall show now that the conclusion of Theorem~\ref{Theorem 4.3} does not 
hold in these cases.   

Assume that $k_\lambda (z)$ is a reproducing kernel of the type considered in 
\eqref{4.1} and $k_\lambda (\lambda)$ stays bounded as $|\lambda| \sra 1$.  
The 
$\sum_{n=0}^\infty \, a_n < \infty$.  Thus the power series for $k_\lambda (z) 
= \sum_{n=0}^\infty \, a_n \left< z,\lambda \right>^n$ converges absolutely 
and uniformly on $\ol{B}_d \times \ol{B}_d$, and for $\lambda \in \ol{B}_d$, 
$k_\lambda \in \s{H}(k)$.  It follows that all functions in $\s{H}(k)$ extend 
to be continuous on $\ol{B}_d$, where $f(\lambda) = \left< f, k_\lambda 
\right>$ for all $f \in \s{H}(k)$, $\lambda \in \ol{B}_d$.  Let $z \in \prtl 
B_d$ and $\s{M} = \{f \in \s{H}(k) : f(z) = 0 \} = \{k_z\}^\perp$.  Clearly 
this is an invariant subspace.  

Let $\{\varphi_n\}$ be the sequence that is associated with $\s{M}$ according 
to Theorem~\ref{1.1}.  Then for $\lambda \in B_d$, $\sum_{n \geq 1} \, 
|\varphi_n (\lambda)|^2 = \ds{\frac{\|P_\s{M} k_\lambda\|^2}{\|k_\lambda\|^2}} 
= 1 - \ds{\frac{|k_\lambda(z)|^2}{k_\lambda (\lambda) \, k_z(z)}}$.  This is a 
continuous function on $\ol{B_d}$, which is zero at $\lambda = z$.  Thus its 
boundary values cannot be zero a.e. $[\sigma]$.  
\end{exam}

\begin{conj} \label{Conjecture 4.5}
If $\{a_n\}$ and $\{b_n\}$ are related to one another as in \eqref{4.1}, and 
if $\sum_{n=1}^\infty \, b_n = 1$, then $\lim_{n \sra \infty} \, a_n/a_{n+1} = 
1$.  
\end{conj}

Note that for kernels considered in \eqref{4.1} the condition 
$k_\lambda(\lambda) \sra \infty$ as $|\lambda| \sra 1$ is equivalent to 
$\sum_{n=1}^\infty = 1$.  Thus, if the conjecture were true, then the 
hypothesis $\lim_{n \sra \infty} \, a_n/a_{n+1} = 1$ in Theorem~\ref{Theorem 
4.3} would be automatically satisfied, and it could be dropped from the 
statement of the theorem.  In support of Conjecture~\ref{Conjecture 4.5} we 
prove the following proposition.  

\begin{prop} \label{Proposition 4.6}
Suppose $\sum_{n=1}^\infty \, b_n = 1$, and either $\sum_{n=1}^\infty \, nb_n 
< \infty$ or $\{a_n\}$ is eventually nonincreasing.  Then $a_n/a_{n+1} \sra 1$ 
as $n \sra \infty$.  
\end{prop}

\begin{proof}
Suppose first that $\sum_{n=1}^\infty \, b_n = 1$ and $c = \sum_{n=1}^\infty 
\, nb_n < \infty$.  For $z \in \mbb{D}$ let $g(z) = 
\ds{\frac{1-\sum_{n=1}^\infty \, b_n z^n}{1-z}} = \sum_{n=1}^\infty \, b_n 
\frac{1-z^n}{1-z} = \sum_{n=1}^\infty \, b_n \, \sum_{j=0}^{n-1} z_j$.  Then 
$|g(z)| \leq c$, and, in fact, $g$ has a Taylor series that converges 
absolutely and uniformly in $\ol{\mbb{D}}$.  Note that $g(1) = c > 0$ and $\Re 
(1-z) g(z) = 1 - \sum_{n=1}^\infty \, b_n \Re z^n \geq b_1(1- \Re z)$, so 
$g(z) \neq 0$ in $\ol{\mbb{D}}$.  It follows from Wiener's Lemma (see 
\cite{rudinfa}, Theorem 11.6) that the function $h(z) = 
\ds{\frac{1}{g(z)}} - \ds{\frac{1}{c}}$ has a Taylor series $h(z) = 
\sum_{k=0}^\infty \, \hat{h}(k) z^k$ with $\sum_{k=0}^\infty \, |\hat{h}(k)| <
\infty$, and $\sum_{k=0}^\infty \, \hat{h}(k) = h(1) = 0$.  

We have $h(z) + 1/c = \ds{\frac{1}{g(z)}} = (1-z) \sum_{n=0}^\infty \, a_n z^n 
= 1+ \sum_{n=1}^\infty \, (a_n - a_{n-1}) z^n$.  We compare coefficients and 
evaluate at $z=1$ to obtain $1/c = 1 + \sum_{n=1}^\infty \, (a_n - a_{n-1}) = 
\lim_{n \sra \infty} \, a_n$, hence $a_n/a_{n+1} \sra 1$ as $n \sra \infty$.  

Now suppose that $\sum_{n=1}^\infty \, b_n =1$ and that $a_n \geq a_{n+1}$ 
for 
all $n \geq N$.  The since $\sum_{n=0}^\infty \, a_n z^n = 
\ds{\frac{1}{1-\sum_{n=1}^\infty \, b_n z^n}}$ we may multiply through and 
compare coefficients.  We obtain $a_{n+1} = \sum_{k=1}^{n+1} \, b_k 
a_{n+1-k}$.  Thus for $n \geq N$ we have $1 \geq a_{n+1}/a_n = 
\sum_{k=1}^{n+1} \, b_k a_{n+1-k}/a_n \geq \sum_{k=1}^{n-N} \, b_k \sra 1$ as 
$n \sra \infty$. 
\end{proof}

\vspace{1cm}

\section{An application to contractive Hilbert modules}

In \cite{arvcurvinv} Arveson defines a contractive Hilbert module to be a 
Hilbert space $\s{H}$ which is also a module over $\s{A} = \mbb{C}[z_1, 
\ldots, z_d]$, the algebra of complex polynomials in $d$ variables, and has 
the property that 

\begin{equation} \label{5.1}
\|z_1 \xi_1 + \cdots + z_d \xi_d\|^2 \leq \|\xi_1\|^2 + \cdots + \|\xi_d\|^2 
\; \text{~for all~} \; \xi_1, \ldots, \xi_d \in \s{H}.  
\end{equation}

Hence the actions of $z_1, \ldots, z_d$ define bounded linear operators on 
$\s{H}$ which we denote by $T_1, \ldots, T_d$, respectively.  From the 
definition of a contractive Hilbert module and \eqref{5.1} the $d$-tuple of 
operators $(T_1, \ldots, T_d)$ satisfies the properties

\begin{equation*}
\sum_{i=1}^\infty \, T_i T_i^* \leq 1_\s{H} \; \text{and} \; T_i T_j = T_j T_i 
\; \text{~for all~} \; 1 \leq i,j \leq d.  
\end{equation*}

\noindent Such $d$-tuples have been called $d$-contractions or row 
contractions.  Associated with any Hilbert module $\s{H}$ there is a 
completely positive map $\Psi : \s{B}(\s{H}) \lra \s{B}(\s{H})$ defined by 
$\Psi(X) = \sum_{i=1}^d \, T_i X T_i^*, \, X \in \s{B}(\s{H})$.  A Hilbert 
module is said to be pure if $\lim_{n \sra \infty} \, \Psi^n(1_\s{H}) = 0 
\sot$.  The rank of $\s{H}$ is defined as the rank of the defect operator 
$\Delta = (1_\s{H} - \Psi(1_\s{H}))^{1/2}$.  

Of course, the spaces of the form $\s{H} = H_d^2(\s{D})$ come with a natural 
Hilbert module structure:  If $\xi = f \in H_d^2(\s{D})$, then $z_i \xi = 
M_{z_i} f, \, i = 1, \ldots, d$.  One verifies that $H_d^2 (\s{D})$ is 
contractive, pure, and rank $H_d^2(\s{D}) = \dim \s{D}$.  These modules serve 
a universal role in the categroy of pure contractive Hilbert modules.  The 
following theorem makes this precise.  

\begin{thm} \label{Theorem 5.1}
Let $\s{H}$ be a pure contractive Hilbert module and let $\s{D}$ be a Hilbert 
space with $\dim \s{D} = \rank \s{H}$.  Then there exists a coisometric module 
homomorphism $U: H_d^2(\s{D}) \lra \s{H}$ that is minimal in the sense that 
$U^* \s{H}$ generates $H_d^2(\s{D})$ as a Hilbert module.  Furthermore, if 
$U': H_d^2 (\s{D}') \lra \s{H}$ is another such map, then there exists a 
unitary operator $V: \s{D} \lra \s{D}'$ such that $U = U'\tilde{V}$, where 
$\tilde{V}(fx) = f Vx$ for all $f \in H_d^2, \, x \in \s{D}$.  
\end{thm}

Theorem~\ref{Theorem 5.1} is a well-known result in dilation theory.  For 
example, it can easily be derived from the results in \cite{agler1}, or, for a 
precise statement in the language of Hilbert modules, see \cite{arvmvot}.  In 
fact, if $k$ is a complete NP kernel, then one can define a category of 
Hilbert modules where the spaces of the type $\s{H}(k, \s{D})$ play the role 
of the universal object.  In this case one uses the functions $\{b_n\}$ of 
\eqref{1.1} to define the completely positive map $\Psi$, and an analogue of 
the above theorem holds (see \cite{agler1}).  

Thus, any pure contractive Hilbert module $\s{H}$ can be associated with a 
submodule $\s{M} = \ker U$ of $H_d^2(\s{D})$.  It follows from 
Lemma~\ref{Lemma 4.1} that $\s{M}$ is a multiplier invariant subspace of 
$H_d^2(\s{D})$, so by Theorem~\ref{Theorem 1.3} there exists an auxiliary 
Hilbert space $\s{E}$ and an inner multiplication operator $\Phi \in 
\s{B}(H_d^2(\s{E}), H_d^2(\s{D}))$ with associated operator vauled multiplier 
$\phi \in \s{B}(\s{E}, \s{D})$ such that $\s{M} = \Phi H_d^2(\s{E})$.  

The {\it curvature invariant} of a finite rank Hilbert module $\s{H}$ was 
introduced in \cite{arvcurvinv}.  To review the definition we need to fix some 
more notation.  If $T_1, \ldots, T_d$ are the operators associated with 
$\s{H}$, then for $\lambda \in B_d$ we set $T(\lambda) = \ol{\lambda}_1 T_1 + 
\cdots \ol{\lambda}_d T_d$.  Since $\s{H}$ has finite rank the space $\Delta 
\s{H}$ is finite dimensional.  We define a $\s{B}(\Delta \s{H})$-valued 
function on $B_d$ by 

\begin{equation*} 
F(\lambda) = (1-|\lambda|^2) \Delta (1_\s{H} - T(\lambda)^*)^{-1} (1_\s{H} - 
T(\lambda))^{-1} \Delta.
\end{equation*}

It can be shown that $F(\lambda)$ is unitarily equivalent to $1_\s{D} - 
\phi(\lambda) \phi(\lambda)^*$, where $\phi$ is the operator valued multiplier 
as in the previous paragraph (see \cite{arvcurvinv}).  Thus, the radial limit 
(or even $K$-limit) of $F$ exists in the strong operator topology for a.e. $z 
\in \prtl B_d$.  The curvature invariant of $\s{H}$ is defined as 

\begin{equation*}
K(\s{H}) = \int_{\prtl B_d} \, \trace F(z) \, d\sigma(z).  
\end{equation*}

It is clear that $0 \leq K(\s{H}) \leq \rank \s{H}$, and it follows that 
$K(\s{H}) = \int_{\prtl B_d} \, \trace (1_\s{D} - \phi(z) \phi(z)^*) \, 
d\sigma(z)$.  The following theomem resolves Problem 1 of \cite{arvcurv}.  

\begin{thm} \label{Theorem 5.2}
If $\s{H}$ is a contractive, pure Hilbert module of finite rank, then 
$K(\s{H})$ is an integer.  

In particular, if $\s{M}$ is the multiplier invariant subspace associated with 
$\s{H}$ as above, then 

\begin{multline*}
K(\s{H}) = \rank \s{H} - \sup \{\dim E_\lambda \s{M} : \lambda \in B_d \} \\
= 
\inf \{\dim \s{D}_\lambda \cap \s{M}^\perp : \lambda \in B_d \}, 
\end{multline*}

\noindent where $\s{D}_\lambda = k_\lambda \s{D}$.  
\end{thm}

\begin{proof}
This follows immediately from Theorem~\ref{Theorem 4.3}.  Recall that for 
$\lambda \in B_d, \, E_\lambda : H_d^2(\s{D}) \lra \s{D}$ denotes the point 
evaluation, $E_\lambda f = f(\lambda)$.  

For a.e. $z \in \prtl B_d, \, F(z) = 1_\s{D} - \phi(z) \phi(z)^*$ is a 
projection of rank $\rank \s{H} - \sup \{\dim E_\lambda \s{M} : \lambda \in 
B_d\} = \inf \{\dim (E_\lambda \s{M})^\perp : \lambda \in B_d\} = \inf \{\dim 
\s{D}_\lambda \cap \s{M}^\perp : \lambda \in B_d \}$, since one easily sees 
that $(E_\lambda \s{M})^\perp = \s{D}_\lambda \cap \s{M}^\perp$.  Hence 
$K(\s{H}) = \int_{\prtl B_d} \, \trace F(z) \, d\sigma(z) = \inf \{\dim 
\s{D}_\lambda \cap \s{M}^\perp : \lambda \in B_d \}$. 
\end{proof}

\begin{exam} \label{Example 5.3}
Let $\s{H}$ be a pure contractive Hilbert module of finite rank, let $U$ be as 
in Theorem~\ref{Theorem 5.1}, and, as above, set $\s{M} = \ker U$.  

\begin{enumerate}
\item \label{aa}
If $\varphi$ is a nonzero scalar multiplier of $H_d^2$ such that $\varphi 
\s{H} = 0$, then $K(\s{H}) = 0$.  
\item \label{bb}
If $\s{M}$ is generated by a family of functions $\{f_n\}_{n \geq 1}$ such 
that there is a nonempty open set $\Omega \appa B_d$ such that the dimension 
of the linear span of $\{f_n(\lambda)\}_{n \geq 1}$ in $\s{D}$ equals $m$ for 
each $\lambda \in \Omega$, then $K(\s{H}) = \rank \s{H} - m$.  
\end{enumerate}
\end{exam}

\begin{proof}
\begin{enumerate} 
\item 
Let $\lambda \in B_d$ such that $\varphi(\lambda) \neq 0$.  
According to 
Theorem~\ref{Theorem 5.2} it suffices to show that $\s{D}_\lambda \cap 
\s{M}^\perp = (0)$.  Thus, let $x \in \s{D}$ be such that $k_\lambda x \in 
\s{M}^\perp$.  Then for any $y \in \s{D}, \, \varphi y$ is a function in 
$\s{M}$ since $\varphi \s{H} = 0$.  Hence $0 = \left< \varphi y, k_\lambda x 
\right> = \varphi(\lambda) \left< y, x \right>_\s{D}$ and it follows that 
$x= 0$.
\item It follows from the hypothesis that $\dim E_\lambda \s{M} = m$ for each 
$\lambda \in \Omega$.  Thus, it follows from Lemma~\ref{Lemma 3.1} that $\sup 
\{\dim E_\lambda \s{M} : \lambda \in B_d \} = m$, and the result follows from 
Theorem~\ref{Theorem 5.2}.  
\end{enumerate}
\end{proof}

It is sometimes possible for $K(\s{H})$ to be defined and finite when $\dim 
\, \s{D} = \infty$.  Thus, a more general resolution of Problem 1 in 
\cite{arvcurv} would follow from an answer to 

\begin{question}
Is $\dim (\ran \phi(z))^\perp$ almost everywhere equal to a constant even if 
$\dim \s{D} = \infty$?  
\end{question}

Note that Theorem~\ref{Theorem 4.3} implies that $\dim \ran \phi(z)$ is a.e. 
equal to a constant.  

%\bibliography{mybib}
\providecommand{\bysame}{\leavevmode\hbox to3em{\hrulefill}\thinspace}

\end{document}